\documentclass[12pt,a4paper]{article}
\usepackage{latexsym,amssymb,amsmath,amsfonts}
\newcommand{\bea}{\begin{eqnarray}}
\newcommand{\eea}{\end{eqnarray}}
\newcommand{\bean}{\begin{eqnarray*}}
\newcommand{\eean}{\end{eqnarray*}}
\newcommand{\nbea}{\begin{neweqnarray}} 
\newcommand{\neea}{\end{neweqnarray}} 

 \newcommand{\qed}{ \mbox
{$\, \Box \,$}} \newenvironment{proof}{\noindent {\bf
Proof.}}{\qed}
\newcounter{secnum}
\newcommand{\sect}[1]
{\protect\section{#1} \protect\setcounter{secnum}{\value{section}}
\protect\setcounter{equation}{0}
\protect\renewcommand{\theequation}{\mbox{\arabic{secnum}.\arabic{equation}}
}} 
\newtheorem{theorem}{Theorem}[section]

\newtheorem{lemma}{Lemma}[section]

\newtheorem{remark}{Remark}[section]
\newcommand{\beq}{\begin{equation}}
\newcommand{\eeq}{\end{equation}}
\def\bea{\begin{eqnarray}}
\def\eea{\end{eqnarray}}

\newcommand\F{\mbox{I\kern-2pt F}}

\newcommand\cG{{\cal G}}
\newcommand\cF{{\cal F}}

\newcommand\cP{{\cal P}}

\newcommand\cD{{\cal D}}

\def\text#1{\hbox{#1}}

\def\endproof{\mbox{\ $\qed$}}

\def\E{{\bf E}}
\def\P{{\bf P}}

\def\build #1_#2{\mathrel{\mathop{\kern 0pt #1}\limits_{#2}}}
\newcommand{\zs}[1]{{\mathchoice{#1}{#1}{\lower.25ex\hbox{$\scriptstyle#1$}}
{\lower0.25ex\hbox{$\scriptscriptstyle#1$}}}}

\renewcommand{\theequation}{\thesection.\arabic{equation}}

\def\endproof{\mbox{\ $\qed$}}
\begin{document}
\title{
On asymptotic normality of sequential LS-estimates of unstable autoregressive processes.
  \thanks{This research was carried out at the Department of
Mathematics, the Strasbourg University, France. The second author
is partially supported by the RFFI-DFG-Grant 02-01-04001. }}
\author{{  L. Galtchouk} \\ 
Strasbourg University \\ 
and \\
 V.Konev\\
Tomsk University\\ 
}
\date{}
\maketitle

\begin{abstract} For estimating the unknown parameters in an
unstable autoregressive AR(p), the paper proposes  sequential 
least squares estimates with a special stopping time defined by
the trace of the observed Fisher information matrix. The limiting
distribution of the sequential LSE is shown to be normal for the
parameter vector lying both inside the stability region and on some
part of its boundary in contrast to the ordinary LSE. The asymptotic
normality of the sequential LSE is provided by a new property of the
observed Fisher information matrix which holds both inside the 
stability region of AR(p) process and on the part of its boundary.
The asymptotic distribution of the stopping time is derived.
 \end{abstract}



\bibliographystyle{plain}

\newpage

\sect{Introduction}

Consider the autoregressive $AR(p)$ model
\begin{equation} \label{1.1}
x_n=\theta_1 x_{n-1}+\ldots +\theta_p x_{n-p}+\varepsilon_n,\,
n=1,2,\ldots,
\end{equation}
where $(x_n)$ is the observation, $(\varepsilon_n)$ is the noise which is a sequence
of independent
 identically distributed (i.i.d.) random variables with $\E\varepsilon_1 = 0$
 and  $0<\E\varepsilon_1^2=\sigma^2<\infty$, $\sigma^2$ is known, $x_0=x_{-1}=\ldots=x_{1-p}=0$; 
 parameters of the model $\theta_1, \ldots,\theta_p$
 are unknown. This model can be expressed in vector form as
 \begin{equation}\label{1.2}
 X_n\,=\,AX_{n-1}\,+\,\xi_n\,,
 \end{equation}
 where $X_n=(x_n,x_{n-1},\ldots,x_{n-p+1})^{\prime}, \xi=(\varepsilon_n,0,\ldots,0)^{\prime}$,
 \begin{equation}\label{1.3}
 A\,=\, 
 \left(
 \begin{array}{cc}
 \theta_1\ldots  & \theta_p \\
 I_{p-1}         & 0
 \end{array}
 \right),
 \end{equation}
 the prime denotes the transposition.

 A commonly used estimate of the parameter vector $\theta=(\theta_1,\ldots,\theta_p)^{\prime}$
 is the least squares estimate (LSE)
\begin{equation}\label{1.4}
\theta(n) =M_n^{-1}\sum_{k=1}^n
X_{k-1 }x_k,\quad
 M_n=\sum_{k=1}^n X_{k-1}X_{k-1}^{\prime},
\end{equation}
where  
 $M_n^{-1}$ denotes the inverse of matrix $M_n$ if $\det M_n>0$ and
$M_n^{-1}=0$ otherwise. Let
\begin{equation}\label{1.5}
 \mathcal{P}(z)=z^p-\theta_1z^{p-1}-\ldots-\theta_p
\end{equation} 
denote the characteristic polynomial of the autoregressive model (\ref{1.1}).
The process (\ref{1.1}) is said to be stable (asymptotically stationary) if
all roots $z_i=z_i(\theta)$ of the characteristic polynomial (\ref{1.5})
lie inside the unit circle, that is the parameter vector $\theta=(\theta_1,\ldots,\theta_p)^{\prime}$ belongs to the parametric stability
region $\Lambda_p$ defined as
\begin{equation}\label{1.6}
 \Lambda_p=\{\theta\in R^p:|z_i(\theta)|<1, i=1,\ldots,p\}\,.
\end{equation} 
The process (\ref{1.1}) is called unstable if the roots of $\mathcal{P}(z)$
lie on or inside the unit circle, that is, $\theta\in[\Lambda_p]$; $[\Lambda_p]$
denotes the closure of the stability region $\Lambda_p$.

It is well known (see,e.g. Anderson (1971), Th.5.5.7) that the LSE $\theta(n)$
is asymptotically normal for all $\theta\in\Lambda_p$, that is
$$
\sqrt{n}(\theta(n)-\theta)\stackrel{\mathcal{L}}{\Longrightarrow}\mathcal{N}(0,F),\
  \hbox{\rm as}\ n\to\infty,
$$ 
where $F=F(\theta)$ is a positive definite matrix, $\stackrel{\mathcal{L}}{\Longrightarrow}$
indicates convergence in law. It should be noted that the asymptotic normality of
$\theta(n)$ is provided by the following asymptotic property of the observed Fisher
information matrix
\begin{equation}\label{1.7}
 \lim_{n\to\infty}M_n/n\,=\,F\\\ \hbox{\rm a.s.}
\end{equation}
for all $\theta\in\Lambda_p$. On the boundary $\partial\Lambda_p$ of the stability
region $\Lambda_p$, this property does not hold and the distribution of $\theta(n)$
is no longer asymptotically normal. The investigation of the asymptotic distribution of LSE
$\theta(n)$ when $x_n$ is unstable goes back to the late fifties with the paper of White (1958)
(see also Ahtola and Tiao (1987), Dickey and Fuller (1979), Rao (1978), Sriram (1987),(1988))
who considered the AR(1) model with i.i.d. $\mathcal{N}(0,\sigma^2)$ random errors
$\varepsilon_n$ and $\theta_1=1$ and established that 
$$
n(\theta(n)-1)\stackrel{\mathcal{L}}{\Longrightarrow}(W^2(1)-1)/\int_0^1\,W^2(t)dt\,,
$$
where $W(t)$ is a standard brownian motion. Subsequently the research of the limiting
distribution of $\theta(n)$ for unstable AR(p) processes has been receiving considerable
attention due to important applications in time series analysis, in  modeling economic and financial data and in system identification and control. We can not go into the detail here
and refer the reader to the paper by Chan and Wei (1988) who derived the limiting
distribution of LSE $\theta(n)$ for the general unstable AR(p) model. By making use of the
functional central limit theorem approach, Chan and Wei expressed the limiting
distribution of LSE $\theta(n)$ in terms of functionals of standard brownian motions. However,
the closed forms of the distribution functions of these functionals are not known and that may
cause difficulties in practice (see section 4 in Chan and Wei).

For the unstable AR(1) model with i.i.d. random errors and $-1\le\theta_1\le 1$, Lai and Siegmund
(1983) proposed to use the sequential least squares estimate for $\theta_1$ which is obtained
from the LSE
$$
\theta_1(n)\,=\,\sum_{k=1}^n\,x_{k-1}x_k/\sum_{k=1}^n\,x_{k-1}^2
$$
by replacing $n$ with a special stopping time $\tau$ based on the observed Fisher information.
They proved that, in contrast with the ordinary LSE $\theta_1(n)$, the sequential LSE is 
asymptotically normal uniformly in $\theta\in[-1,1]$. For the unstable AR(2) model, Galtchouk 
and Konev
(2006) applied the sequential LSE with a particular stopping time and established that it is
asymptotically normal not only inside the stability rigion $\Lambda$ but also for its boundary
points $\theta$ corresponding to a pair of conjugate complex roots $z_1=e^{i\phi}, z_2=e^{-i\phi}$ of the polynomial (\ref{1.5}).

In this paper, for the case of unstable AR(p) process, we propose a sequential LSE for $\theta$ 
and find the conditions on $\theta$ (see Conditions 1-3 in the next section) ensuring its 
asymptotic normality. The set ${\tilde{\Lambda}}_p$ of the points $\theta$ , satisfying these conditions
 includes the stability region $\Lambda_p$ and some part of its boundary. It is shown that
the convergence of the sequential LSE to the normal distribution is uniform in $\theta\in K$ for
any compact set $K\in{\tilde{\Lambda}}_p$ (see Theorem 2.1). The extension of the property of asymptotic normality
of the sequential estimate to the part of the boundary $\partial\Lambda_p$ is achieved by making
use of a new property of observed Fisher information matrix $M_n$, which holds in a broader subset of 
$[\Lambda_p]$ as compared with (\ref{1.7}) (see Lemma 3.3). 

The remainder of this paper is arranged as follows. In Section 2 we introduce a sequential procedure 
for estimating the parameter vector $\theta=(\theta_1,\ldots,\theta_p)^{\prime}$ in
(\ref{1.1}) and study its properties. Section 3 gives a new property of the observed Fisher 
information matrix and establishes some technical results. In Section 4 we prove Theorem 2.2 from
Section 2 on asymptotic distribution of the stopping time.

\sect{ Sequential least squares estimate. \\ Uniform asymptotic
normality.}

In this section we consider the sequential least squares estimate and study its asymptotic
properties. We define the sequential LSE for the parameter vector 
$\theta=(\theta_1,\ldots,\theta_p)^{\prime}$ in model (\ref{1.1}) as
\begin{equation}\label{2.1}
\theta(\tau(h))=M^{-1}_{\tau(h)}\sum^{\tau(h)}_{k=1}X_{k-1}x_k,
\end{equation}
where
\begin{equation}\label{2.2}
\tau(h)=\inf\left\{n\ge 1: tr\,M_n \ge
h\sigma^2\right\},\ \inf\{\emptyset\}=+\infty,
\end{equation} 
is stopping time, $h$ is a positive number (threshold).

Assume that the parameter vector $\theta=(\theta_1,\ldots,\theta_p)^{\prime}$ in (\ref{1.1})
satisfies the following Conditions.

{\bf Condition 1}. Parameter $\theta=(\theta_1,\ldots,\theta_p)^{\prime}$ is such that all
roots $z_i=z_i(\theta)$ of the characteristic polynomial  (\ref{1.5}) lie inside or on the unite circle.

{\bf Condition 2}. All the roots $z_i=z_i(\theta)$ of $\mathcal{P}(z)$, which are equal to one
in modulus, are simple.

{\bf Condition 3}. The system of linear equations with respect to $Y_1,\ldots,Y_{p-1}$
\begin{equation}\label{2.3}
\left\{
\begin{array}{l}
Y_1-\sum_{l=2}^p\,\theta_lY_{l-1}=\theta_1\\
-\sum_{k=1}^{j-1}\,\theta_{j-k}+Y_j-\sum_{k=1}^{p-j}\,\theta_{k+j}Y_k=\theta_j,\\ 2\le j\le p-1,
\end{array}
\right.
\end{equation}
has a unique solution $(Y_1,\ldots, Y_{p-1}), Y_i=\kappa_i(\theta), 1\le i\le p-1$, and the matrix
\begin{equation}\label{2.4}
L(\theta_1,\ldots,\theta_p)=
\left(
\begin{array}{ccccc}
1                & \kappa_1(\theta) & \kappa_2(\theta) & \ldots & \kappa_{p-1}(\theta) \\
\kappa_1(\theta) & 1                &\kappa_1(\theta)  & \ldots & \kappa_{p-2}(\theta) \\
\vdots           & \vdots           & \vdots           & \ddots & \vdots              \\
\kappa_{p-1}(\theta)&\kappa_{p-2}(\theta)& \ldots      & \kappa_{1}(\theta) & 1       \\
\end{array}
\right)
\end{equation}
is positive definite.

Let ${\stackrel{\circ}{\Lambda}}_p$ denote all $\theta=(\theta_1,\ldots,\theta_p)^{\prime}$ in 
(\ref{1.1}) which satisfy Conditions 1,2, and ${\tilde{\Lambda}}_p$-- all
$\theta=(\theta_1,\ldots,\theta_p)^{\prime}$ satisfying all Conditions 1-3.

{\it Example 2.1.} For AR(2) process, one finds
$$
\Lambda_2=\{\theta=(\theta_1,\theta_2)^{\prime}: -1+\theta_2<\theta_1<1-\theta_2, |\theta_2|<1\}\,,
$$
$$
[\Lambda_2]=\{\theta=(\theta_1,\theta_2)^{\prime}: -1+\theta_2\le\theta_1\le 1-\theta_2, |\theta_2|\le 1\}\,,
$$
$$
{\stackrel{\circ}{\Lambda}}_p=[\Lambda_2]\setminus\{(-2,-1),(2,-1)\}\,,
$$
$$
{\tilde{\Lambda}}_2=\{\theta=(\theta_1,\theta_2)^{\prime}: -1+\theta_2<\theta_1<1-\theta_2, -1\le\theta_2<1\}\,,
$$
$$
L(\theta_1,\theta_2)=
\left(
\begin{array}{cc}
1                & \theta_1/(1-\theta_2)  \\
\theta_1/(1-\theta_2)& 1                \\
\end{array}
\right)\,.
$$
{\it Example 2.2.} By numerical calculation for AR(3) process, one can check that 
Conditions 1-3 are satisfied, for example, for the values of $\theta=(\theta_1,\theta_2,\theta_3)$
 such that $z_1(\theta)=e^{i\phi},\ z_2(\theta)=e^{-i\phi}$ with $3\pi/10\le\phi\le 3\pi/5$ and
 $-1\le z_3(\theta)\le -0.5$.

\vspace{0.5cm}
As is shown in Lemma 3.3 (Section 3), Conditions 1-3, imposed on the parameter
 $\theta=(\theta_1,\ldots,\theta_p)^{\prime}$ in (\ref{1.1}), provide the convergence of the ratio 
 $M_n/\sum_{k=1}^n\,x_{k-1}^2$ to the matrix $L(\theta_1,\ldots,\theta_p)$ given in (\ref{2.4}).
This property can be viewed as an extension of (\ref{1.7}) outside the stability region (\ref{1.6}).

\begin{remark}
It will be observed that $\Lambda_p\subset{\tilde{\Lambda}}_p$ and, for all $\theta\in{\tilde{\Lambda}}_p$, one has
\begin{equation}\label{2.5}
\lim_{n\to\infty}\frac{M_n}{\sum_{k=1}^n\,x_{k-1}^2}=\frac{pF}{tr\,F}
=\Lambda(\theta_1,\ldots,\theta_p)\\\  \hbox{\rm a.s.},
\end{equation} 
where $F$ is the same as in (\ref{1.7}). Indeed, by making use of the identity
\begin{equation}\label{2.6}
\sum_{k=1}^n\,x_{k-1}^2=\frac{1}{p}\sum_{k=1}^n\,\|X_{k-1}\|^2+\frac{1}{p}\sum_{i=2}^p\sum_{l=n-i+2}^n\,x_{l-1}^2\,,
\end{equation} 
one obtains
$$
\frac{M_n}{\sum_{k=1}^n\,x_{k-1}^2}=\frac{M_n}{n}\left(\frac{1}{p}tr\frac{M_n}{n}
(1+(\sum_{k=1}^n\|X_{k-1}\|^2)^{-1}\sum_{i=2}^p\sum_{l=n-i+2}^nx_{l-1}^2)\right)^{-1}.
$$
Limiting $n\to\infty$, one comes to (\ref{2.5}), in view of (\ref{1.7}).
\end{remark}
\begin{theorem}
Suppose that in the AR(p) model (\ref{1.1}), $(\varepsilon_n)$ is a sequence of i.i.d. random
variables with $\E\varepsilon_n=0$ and $\E\varepsilon_n^2=\sigma^2<\infty$ and the parameter vector 
$\theta=(\theta_1,\ldots,\theta_p)^{\prime}$ satisfies Conditions 1-3. Then for any compact set $K\subset{\tilde{\Lambda}}_p$
\begin{equation}\label{2.7}
\lim_{h\to\infty}\sup_{\theta\in K}\sup_{t\in R^p}
\left|\P_{\theta}\left(M_{\tau(h)}^{1/2}(\theta(\tau(h))-\theta)\le t\right)-\Phi_p(\frac{t}{\sigma})\right|=0\,,
\end{equation}
where $\Phi_p(t)=\Phi(t_1)\cdots\Phi(t_p), \Phi$ is the standard normal distribution function,
${\tilde{\Lambda}}_p$ is defined in Condition 3.
\end{theorem}

  {\sl Proof.} Substituting (\ref{1.1}) in (\ref{2.1}) yields
\begin{equation}\label{2.8}
M^{1/2}_{\tau(h)}(\theta(\tau(h)-\theta)=M^{-1/2}_{\tau(h)}\sum^{\tau(h)}_{k=1}X_{k-1}\varepsilon_k
=\sqrt{h}M^{-1/2}_{\tau(h)}L^{1/2}(\theta_1,\ldots,\theta_p)Y_h,
\end{equation}
where 
\begin{equation}\label{2.9}
Y_h=\frac{1}{\sqrt{h}}\sum^{\tau(h)}_{k=1}L^{-1/2}(\theta_1,\ldots,\theta_p)X_{k-1}\varepsilon_k
\end{equation}
and $L(\theta_1,\ldots,\theta_p)$ is given in (\ref{2.4}). Denote
$$
{\cG}_{\tau(h)}=L^{-1/2}(\theta_1,\ldots,\theta_p)M_{\tau(h)}L^{-1/2}(\theta_1,\ldots,\theta_p).
$$
One can easily verify that
$$
\left\|L^{-1/2}(\theta_1,\ldots,\theta_p)\frac{M_{\tau(h)}^{1/2}}{\sqrt{h}}-I_p\right\|^2=
\left\|\frac{1}{\sqrt{h}}{\cG}_{\tau(h)}^{1/2}-I_p\right\|^2
$$
$$
\le\|h^{-1}{\cG}_{\tau(h)}-I_p\|^2\le tr\,L^{-1}(\theta_1,\ldots,\theta_p)
\|h^{-1}M_{\tau(h)}-L(\theta_1,\ldots,\theta_p)\|^2\,.
$$
From here, by making use of Lemma 3.4 from Section 3, one gets, for any compact set $K\subset{\tilde{\Lambda}}_p$
and $\delta>0$,
\begin{equation}\label{2.10}
\lim_{h\to\infty}\sup_{\theta\in K}\P_{\theta}\left(\|\sqrt{h}M_{\tau(h)}^{-1/2}L^{1/2}(\theta_1,\ldots,\theta_p)-I_p\|>\delta\right)=0\,.
\end{equation}
Now we prove that for any compact set $K\subset{\tilde{\Lambda}}_p$ and for each constant vector
$v\in R^p$ with $\|v\|=1$
\begin{equation}\label{2.11}
\lim_{h\to\infty}\sup_{\theta\in K}\sup_{t\in R}|\P_{\theta}(v^{\prime}Y_h\le t)-\Phi(t)|=0\,.
\end{equation}
In view of (\ref{2.9}), one has
$$
v^{\prime}Y_h=\frac{1}{\sqrt{h}}\sum_{k=1}^{\tau(h)}\,g_{k-1}\varepsilon_k,\\\ 
g_{k-1}=v^{\prime}L^{-1/2}(\theta_1,\ldots,\theta_p)X_{k-1}.
$$
For each $h>0$, we define an auxiliary stopping time as
$$
\tau_0=\tau_0(h)=\inf\{n\ge 1: \sum_{k=1}^{n}\,g_{k-1}^2\ge h\},\\ \inf\{\emptyset\}=+\infty.
$$
Further we make use of the representation
$$
v^{\prime}Y_h=\frac{1}{\sqrt{h}}\sum_{k=1}^{\tau_0(h)}\,g_{k-1}\varepsilon_k
+\eta(h)+\Delta(h)\,,
$$
where $\Delta(h)=\sum_{i=1}^4\Delta_i(h)$,
$$
\Delta_1(h)=h^{-1/2}I_{(\tau(h)=1)}g_0\varepsilon_1,
\Delta_2(h)=h^{-1/2}g_{\tau(h)-1}\varepsilon_{\tau(h)}\,,
$$
$$
\Delta_3(h)=-h^{-1/2}I_{(\tau_0(h)=1)}g_0\varepsilon_1,
\Delta_4(h)=h^{-1/2}g_{\tau_0(h)-1}\varepsilon_{\tau_0(h)}\,,
$$
$$
\eta(h)=\frac{1}{\sqrt{h}}\sum_{k=1}^{\tau(h)-1}\,g_{k-1}\varepsilon_k-
\frac{1}{\sqrt{h}}\sum_{k=1}^{\tau_0(h)-1}\,I_{(\tau_0(h)>1)}g_{k-1}\varepsilon_k\,.
$$
Now we show that
\begin{equation}\label{2.12}
\lim_{h\to\infty}\sup_{\theta\in K}\sup_{t\in R}|\P_{\theta}(\frac{1}{\sqrt{h}}\sum_{k=1}^{\tau_0(h)}\,g_{k-1}\varepsilon_k\le t)-\Phi(t)|=0
\end{equation}
and for any $\delta>0$
\begin{equation}\label{2.13}
\lim_{h\to\infty}\sup_{\theta\in K}\P_{\theta}(|\eta(h)|>\delta)=0\,,
\end{equation}
\begin{equation}\label{2.14}
\lim_{h\to\infty}\sup_{\theta\in K}\P_{\theta}(|\Delta(h)|>\delta)=0\,.
\end{equation}
The proof of (\ref{2.12}) is based on Proposition 3.1 from the paper by Lai and Siegmund (1983).
Actually one needs to check only the condition $A_6$, that is, for each $\delta>0$,
\begin{equation}\label{2.15}
\lim_{m\to\infty}\sup_{\theta\in K}\P_{\theta}(g_n^2\ge\delta\sum_{k=1}^n\,g_{k-1}^2
\\\ \hbox{\rm for some}\ n\ge m)=0\,.
\end{equation}
Conditions $A_1-A_5$ are evidently satisfied.
It will be noted that
$$
\sum_{k=1}^n\,g_{k-1}^2=\left(v^{\prime}L^{-1/2}(M_n/\sum_{k=1}^n\,x_{k-1}^2-L)L^{-1/2}v+1\right)\sum_{k=1}^n\,x_{k-1}^2.
$$
Proceeding from this equality one gets the inclusion
$$
\{g_n^2\ge\delta\sum_{k=1}^n\,g_{k-1}^2\\\ \hbox{\rm for some}\ n\ge m\}\subseteq
$$
$$
\subseteq\{\|X_n\|^2\ge\delta_1\sum_{k=1}^n\,g_{k-1}^2\\\ \hbox{\rm for some}\ n\ge m\}
$$
$$
=\{\|X_n\|^2\ge\delta_1\sum_{k=1}^n\,x_{k-1}^2[1+v^{\prime}L^{-1/2}(M_n/\sum_{k=1}^n\,x_{k-1}^2-L)L^{-1/2}v]\\\ \hbox{\rm for some}\ n\ge m\}
$$
$$
\subseteq\{\|X_n\|^2\ge\delta_1\sum_{k=1}^n\,x_{k-1}^2[1-\|L^{1/2}v\|^2\|M_n/\sum_{k=1}^n\,x_{k-1}^2-L\|]\\\ \hbox{\rm for some}\ n\ge m\}
$$
$$
\subseteq\{\|X_n\|^2\ge\delta_1\sum_{k=1}^n\,x_{k-1}^2[1-a^*\|M_n/\sum_{k=1}^n\,x_{k-1}^2-L\|]\\\ \hbox{\rm for some}\ n\ge m\}
$$
$$
\subset\{\|M_n/\sum_{k=1}^n\,x_{k-1}^2-L\|\ge(2a^*)^{-1}\\\ \hbox{\rm for some}\ n\ge m\}
$$
$$
\cup\{\|X_n\|^2\ge\frac{\delta_1}{2}\sum_{k=1}^n\,x_{k-1}^2\\\ \hbox{\rm for some}\ n\ge m\}\,,
$$
where $\delta_1=\delta/a^*, a^*=\sup_{\theta\in K}\|v^{\prime}L^{-1/2}\|^2$.

This, in view of Lemmas 3.1,3.3, yields (\ref{2.15}). It will be observed that (\ref{2.15})
enables one to show (by the same argument as in Lemma 3.5) that, for any compact set $K\subset{\tilde{\Lambda}}_p$ and $\delta>0$
\begin{equation}\label{2.16}
\lim_{h\to\infty}\sup_{\theta\in K}\P_{\theta}(g_{\tau_0-1}^2/\sum_{k=1}^{\tau_0-1}\,g_{k-1}^2\ge\delta)=0\,.
\end{equation}
Now we check (\ref{2.13}). One can easily verify that
$$
\E_{\theta}\eta^2(h)=\E_{\theta}u(h),
u(h)=\frac{1}{h}\left|\sum_{k=1}^{\tau-1}\,g_{k-1}^2-\sum_{k=1}^{\tau_0-1}\,g_{k-1}^2\right|.
$$
The random variable $u(h)$ is uniformly bounded from above uniformly in $\theta\in K$ because
$$
u(h)\le\frac{1}{h}\sum_{k=1}^{\tau-1}\,g_{k-1}^2+1=
\frac{1}{h}v^{\prime}L^{-1/2}M_{\tau(h)-1}L^{-1/2}v+1
$$
$$
\le\frac{a^*}{h}\sum_{k=1}^{\tau-1}\,\|X_{k-1}\|^2+1\le a^*+1\,.
$$ 
Therefore, it suffices to show that for each $\delta>0$
\begin{equation}\label{2.17}
\lim_{h\to\infty}\sup_{\theta\in K}\P_{\theta}(u(h)\ge\delta)=0\,.
\end{equation}
To this end, one can use the following estimate
$$
u(h)=\frac{1}{h}\left|v^{\prime}L^{-1/2}M_{\tau(h)-1}L^{-1/2}v-\sum_{k=1}^{\tau_0-1}\,g_{k-1}^2\right|
$$
$$
=h^{-1}\sum_{k=1}^{\tau-1}\,x_{k-1}^2v^{\prime}L^{-1/2}\left(M_{\tau(h)-1}/\sum_{k=1}^{\tau-1}\,x_{k-1}^2\,-\,L\right)L^{-1/2}v
$$
$$
+h^{-1}\sum_{k=1}^{\tau-1}\,x_{k-1}^2-h^{-1}\sum_{k=1}^{\tau_0-1}\,g_{k-1}^2
$$
$$
\le a^*\|M_{\tau(h)-1}/\sum_{k=1}^{\tau(h)-1}\,x_{k-1}^2\,-\,L\|
+x_{\tau(h)-1}^2/\sum_{k=1}^{\tau(h)-1}\,x_{k-1}^2+g_{\tau_0(h)-1}^2/\sum_{k=1}^{\tau_0(h)-1}\,g_{k-1}^2\,.
$$
From here, by making use of (\ref{2.16}) and Lemmas 3.4,3.5 one comes to (\ref{2.17}) which 
 implies (\ref{2.13}). By a similar argument, one can check (\ref{2.14}). This completes the proof 
 of (\ref{2.11}). Combining (\ref{2.10}) and (\ref{2.11}) one arrives at (\ref{2.7}).
Hence Theorem 2.1. \endproof

Now we will study the properties of the stopping time $\tau(h)$ defined by (\ref{2.2}).
 Further we need the following
functionals
\begin{equation}\label{2.18}
\begin{array}{l}
J_1(x;t)=\int_0^t\,x^2(s)ds,\\
J_2(x,y;t)=\int_0^t\,(x^2(s)+y^2(s))ds,\\
J_3(x,y;t)=\int_0^t\,(x^2(s)+\mu_1y^2(s))ds,\\
J_4(x,y,z;t)=J_2(x,y;t)+\mu_2J_1(z;t),\\
J_5(x,y,z,u;t)=J_2(x,y;t)+\mu_3J_1(z;t)+\mu_4J_1(u;t),
\end{array}
\end{equation}
where $\mu_i, i=\overline{1,4}$ are defined by  (\ref{4.17}),(\ref{4.27}) and (\ref{4.28}).
For the set ${\tilde{\Lambda}}_p$ of the parameter vector $\theta=(\theta_1,\ldots,\theta_p)
^{\prime}$ satisfying Conditions 1-3, we introduce the following subsets
 belonging to its boundary $\partial{\tilde{\Lambda}}_p$
\begin{equation}\label{2.19}
\begin{array}{l}
\Gamma_1(p)=\{\theta\in\partial{\tilde{\Lambda}}_p: z_1(\theta)=-1,
 | z_k(\theta)|<1, k=\overline{2,p}\}\\
\Gamma_2(p)=\{\theta\in\partial{\tilde{\Lambda}}_p: z_1(\theta)=1,
 | z_k(\theta)|<1, k=\overline{2,p}\}\\ 
 \Gamma_3(p)=\{\theta\in\partial{\tilde{\Lambda}}_p: z_1(\theta)=e^{i\phi}, 
 z_2(\theta)=e^{-i\phi}, \phi\in(0,\pi), | z_k(\theta)|<1, k=\overline{3,p}\}\\
 \Gamma_4(p)=\{\theta\in\partial{\tilde{\Lambda}}_p: z_1(\theta)=-1, 
 z_2(\theta)=1,   | z_k(\theta)|<1, k=\overline{3,p}\}\\
 \Gamma_5(p)=\{\theta\in\partial{\tilde{\Lambda}}_p: z_1(\theta)=-1, z_2(\theta)=e^{i\phi}, 
 z_3(\theta)=e^{-i\phi}, \phi\in(0,\pi),\\
  | z_k(\theta)|<1, k=\overline{4,p}\},\\  
 \Gamma_6(p)=\{\theta\in\partial{\tilde{\Lambda}}_p: z_1(\theta)=1, z_2(\theta)=e^{i\phi}, 
 z_3(\theta)=e^{-i\phi}, \phi\in(0,\pi),\\
  | z_k(\theta)|<1, k=\overline{4,p}\}\\
 \Gamma_7(p)=\{\theta\in\partial{\tilde{\Lambda}}_p: z_1(\theta)=-1,  z_2(\theta)=1, z_3(\theta)=e^{i\phi}, 
 z_4(\theta)=e^{-i\phi}, \phi\in(0,\pi),\\
  | z_k(\theta)|<1, k=\overline{5,p}\},       
\end{array}
\end{equation}
where $z_k(\theta)$ are roots of the characteristic polynomial (\ref{1.5}).

It will be noted that all these sets will be used only for the AR(p) model (\ref{1.1})
with $p\ge 5$. In the case when $p\le 4$, it is obious which of the sets $\Gamma_i(p)$ are
odd and how to amend the remaining subsets $\Gamma_i(p)$.
 
\begin{theorem}\label{T.2.2} Suppose that in the AR(p) model (\ref{1.1}), $(\varepsilon_n)_{n\ge 1}$
 is a sequence of i.i.d. random variables with $\E\varepsilon_n=0,\ 0<\E\varepsilon_n^2=\sigma^2<\infty$
 and the parameter vector $\theta=(\theta_1,\ldots,\theta_p)$ satisfies Conditions 1-3. Let   
 $\tau(h)$ be defined by (\ref{2.2}).    Then, for each $\theta\in\Lambda_p$,
\begin{equation}\label{2.20}
   \P_{\theta}-\lim_{h\to\infty}\frac{\tau(h)}{h}=\frac{\sigma^2}{tr F}.
\end{equation}
Moreover, for each $\theta \in\partial{\tilde{\Lambda}}_p$, as $h\to\infty$,
 \begin{equation}\label{2.21}
    \frac{\tau(h)}{b_i\sqrt{h}}
    \stackrel{\mathcal{L}}{\Longrightarrow}\nu_i,\ \hbox{\rm if}\ \theta\in\Gamma_i(p),
    1\le i\le 7\,,
\end{equation}
where $\Lambda_p$ is given in (\ref{1.6});
$$
\begin{array}{l}
\nu_1 = \inf\left\{t\ge 0: J_1(W_1;t)\ge 1\right\},\\ 
\nu_2 = \inf\left\{t\ge 0: J_1(W_2;t)\ge 1\right\}, \\
\nu_3 = \inf\left\{t\ge 0: J_2(W_1,W_2;t)\ge 1\right\},\\ 
\nu_4 = \inf\left\{t\ge 0: J_3(W_1,W_2;t)\ge 1\right\},\\ 
\nu_i = \inf\left\{t\ge 0: J_4(W_1,W_2,W_3;t)\ge 1\right\},\ i=5,6,\\ 
\nu_7 = \inf\left\{t\ge 0: J_5(W_1,W_2,W_3,W_4;t)\ge 1\right\};
\end{array}
$$ 
$b_1,\ldots,b_7$ are defined by (\ref{4.7}),(\ref{4.8}),(\ref{4.10}),(\ref{4.17}),
(\ref{4.27}) and (\ref{4.28}), respectively; $W_1(t),\ldots,W_4(t)$ are independent
standard brownian motions.
\end{theorem}
The proof of Theorem 2.2 is given in the Appendix.

\sect{Auxiliary propositions.}

In this Section we establish some properties of the process
(\ref{1.1}) and the observed Fisher information matrix
$M_n$ used in Section 2.

We need some notations. Let $z_1(\theta),\ldots,z_q(\theta)$ denote all the distinct roots of
the characteristic polynomial (\ref{1.5}), $m_i(\theta)$ be the multiplicity of $z_i(\theta),\ (\sum_{i=1}^qm_i(\theta)=p)$. Let
$$
\rho(\theta)=
\left\{
\begin{array}{cc}
\max(m_i(\theta): |z_i(\theta)|=1)& \hbox{\rm if}\ \max|z_i(\theta)|=1,\\
0 &  \hbox{\rm if}\ \max|z_i(\theta)|\not=1\   \hbox{\rm for all}\ i=\overline{1,q}.
\end{array}
\right.
$$
Formally the set ${\stackrel{\circ}{\Lambda}}_p$ introduced in Condition 3 can be written as
\begin{equation}\label{3.1}
{\stackrel{\circ}{\Lambda}}_p=\Lambda_p\cup\{\theta: \max_{1\le i\le q}|z_i(\theta)|=1,\ \rho(\theta)=1\}.
\end{equation}
It includes both the stability region $\Lambda_p$ and the points $\theta$ of its boundary for which
all the roots of the polynomial (\ref{1.5}), lying on the unit circle, are simple. 
\begin{lemma}\label{l.3.1}
Let $(x_n)_{n\ge 0}$ be an autoregressive process defined by (\ref{1.1}). Then for any compact set $K\subset{\stackrel{\circ}{\Lambda}}_p$ and $\delta>0$
$$
\lim_{m\to\infty}\sup_{\theta\in K}\P_{\theta}(\max_{0\le i\le p-1}x_{n-i}^2\ge\delta\sum_{k=1}^n\,x_{k-1}^2\\\ \hbox{\rm for some}\ n\ge m)=0\,.
$$
\end{lemma}
{\sl Proof.} Taking into account the equation (\ref{1.2}), it suffices to show that for any compact set $K\subset{\stackrel{\circ}{\Lambda}}_p$ and $\delta>0$
$$
\lim_{m\to\infty}\sup_{\theta\in K}\P_{\theta}(B_m(\delta))=0\,,
$$
where
$$
B_m(\delta)=\{\|X_n\|^2\ge\delta\sum_{k=1}^n\,\|X_{k-1}\|^2\\\ \hbox{\rm for some}\ n\ge m\}.
$$
Now we estimate the ratio $\|X_n\|^2/\sum_{k=1}^n\,\|X_{k-1}\|^2$ from above. For each $1\le s<n$, 
we introduce the quantity
$$
l_s=\min\{1\le i\le s: \min_{1\le j\le s}\|X_{k-j}\|^2=\|X_{k-i}\|^2\}
$$
and have the inequality
\begin{equation}\label{3.2}
\sum_{k=1}^n\,\|X_{k-1}\|^2\ge\sum_{k=1}^s\,\|X_{n-k}\|^2\ge s\|X_{n-l_s}\|^2\,.
\end{equation}
On the other hand, it follows from (\ref{1.1}) that
$$
X_n=A^{l_s}X_{n-l_s}+\sum_{i=0}^{l_s-1}A^i\xi_{n-i}
$$
and, therefore, one gets
$$
\|X_n\|^2\le 2\|A^{l_s}\|^2\|X_{n-l_s}\|^2+2\|\sum_{i=0}^{l_s-1}A^i\xi_{n-i}\|^2
$$
$$
\le 2\|A^{l_s}\|^2\|X_{n-l_s}\|^2+2\left(\sum_{i=0}^{s-1}\|A^i\xi_{n-i}\|\right)^2
\le 2\|A^{l_s}\|^2\|X_{n-l_s}\|^2+2s\sum_{i=1}^{s-1}\|A^i\xi_{n-i}\|^2
$$
\begin{equation}\label{3.3}
\le 2\|A^{l_s}\|^2\|X_{n-l_s}\|^2+2s\sum_{i=1}^{s-1}\|A^i\|^2\varepsilon_{n-i}^2\,.
\end{equation}
Further it will be observed that, for every compact set $K\subset{\stackrel{\circ}{\Lambda}}_p$,
there exists a positive number $\kappa$ such that
\begin{equation}\label{3.4}
 \sup_{\theta\in K}\|A^n\|^2\le\kappa,\\\ n\ge 1\,.
\end{equation}
Indeed, we express $A$ in its Jordan normal form
\begin{equation}\label{3.5}
  A=S{\cD}S^{-1}\,,
\end{equation}
where $\cD=\hbox{\rm diag}(J_1,\ldots,J_q), J_l$ is the $m_l\times m_l$ submatrix of the form
$$
J_l=
\left(
\begin{array}{ccccc}
z_l & 1 & 0 &\ldots & 0\\
0 & z_l & 1& \ldots & 0\\
0 &\ldots &  {} & z_l & 1\\
0& \ldots & {} & {} & z_l
\end{array}
\right)
$$
if $z_l$ is a multiple root with multiplicity $m_l\ge 2$, and $J_l=z_l$ if $z_l$ is a simple
root. 

By direct computation with (\ref{3.5}) one finds $A^n=S{\cD}^nS^{-1},\ {\cD}^n=\hbox{\rm diag}(J_1^n,\ldots,J_q^n)$, 
where the powers of the matrix $J_l$ are equal to $z_l^n$ for a simple root $z_l$ and consist 
of the elements (see, R.Varga (2000))
$$
<J_l^n>_{ij}=
\left\{
\begin{array}{cc}
0, & j<i\,, \\
\left(\stackrel{n}{j-i}\right)z_l^{n-j+i}, &  i\le j\le\min(m_l,n+i)\,, \\
0, & n+i<j\le m_l\,, 
\end{array}
\right.
$$
for the roots $z_l$ with multiplicity $m_l\ge 2$, $\left(\stackrel{n}{j-i}\right)$ is the binomial coefficient. From here, in view of the definition (\ref{3.1}), one comes to (\ref{3.4}). By making use of (\ref{3.3}) and (\ref{3.4}), one obtains 
$$
\|X_n\|^2\le 2\kappa \|X_{n-l_s}\|^2+ 2s\kappa\sum_{i=0}^{s-1}\varepsilon_{n-i}^2\,.
$$
Combining this inequality and (\ref{3.2}) yields
$$
\frac{\|X_n\|^2}{\sum_{k=1}^n\|X_{k-1}\|^2}\le\frac{2\kappa}{s}+2s\kappa \frac{\sum_{i=0}^{s-1}\varepsilon_{n-i}^2}{\sum_{k=1}^n\|X_{k-1}\|^2}\,.
$$
It remains to use elementary inequality
$$
\sum_{i=0}^{n-1}\varepsilon_k^2\le 2(1+\|A\|^2)\sum_{k=1}^n\|X_{k-1}\|^2\,,
$$
which follows from (\ref{1.2}), to derive the desired estimate for the ratio
$$
\frac{\|X_n\|^2}{\sum_{k=1}^n\|X_{k-1}\|^2}\le\frac{2\kappa}{s}+2s\kappa \frac{2(1+\|A\|^2)\sum_{i=0}^{s-1}\varepsilon_{n-i}^2}{\sum_{k=1}^{n-1}\varepsilon_k^2}\,.
$$
This inequality implies the inclusion
$$
B_m(\delta)\subset\{2\kappa/s>\delta/2\\\ \hbox{\rm for some}\ n\ge m\}
$$
$$
\cup\{4(1+\|A\|^2)s\kappa\sum_{i=0}^{s-1}\varepsilon_{i-1}^2/\sum_{k=1}^{n-1}\varepsilon_k^2>\frac{\delta}{2}\\\ \hbox{\rm for some}\ n\ge m\}\,.
$$
Therefore, for sufficiently large $s$, one gets
$$
\sup_{\theta\in K}\P_{\theta}(B_m(\delta)\le\P\{\nu\sum_{i=0}^{s-1}\varepsilon_{i-1}^2/\sum_{k=1}^{n-1}\varepsilon_k^2>\frac{\delta}{2}\\\ \hbox{\rm for some}\ n\ge m\}\,,
$$
where $\nu=\sup_{\theta\in K}4(1+\|A\|^2)s\kappa$. Limiting $m\to\infty$ and applying the law of large numbers one comes to the assertion of Lemma 3.1.

\begin{lemma}\label{l.3.2}
Let $(x_n)_{n\ge 0}$ be an autoregressive process defined by (\ref{1.1}). Then for any compact 
set $K\subset{\stackrel{\circ}{\Lambda}}_p$ and each $l=\overline{1,p-1}$
\begin{equation}\label{3.6}
\lim_{m\to\infty}\sup_{\theta\in K}\P_{\theta}(\left|\sum_{k=1}^n\,x_{k-l}\varepsilon_k\right| \ge\delta\sum_{k=1}^n\,x_{k-1}^2\\\ \hbox{\rm for some}\ n\ge m)=0\,,
\end{equation}
where ${\stackrel{\circ}{\Lambda}}_p$ is given in (\ref{3.1})
\end{lemma}
{\sl Proof.} We will apply Lemma 2.2 from the paper by Lai and 
Siegmund (1983)) given in the Appendix.
Let $c_n=n^{3/4}$. For the set of interest one has
the following inclusions 
$$
\left\{|\sum_{k=1}^nx_{k-l}\varepsilon_k|>\delta\sum_{k=1}^nx_{k-1}^2\
\hbox{\rm for some}\ n\ge m\right\} 
$$
$$
\subseteq \left\{|\sum_{k=1}^nx_{k-l}\varepsilon_k|>\delta\sum_{k=1}^nx_{k-l}^2\
\hbox{\rm for some}\ n\ge m\right\} 
 $$
 $$
=\left\{\frac{|\sum_{k=1}^nx_{k-l}\varepsilon_k|}{(\sum_{k=1}^nx_{k-l}^2)^{2/3}\vee
c_n} \cdot\frac{(\sum_{k=1}^nx_{k-l}^2)^{2/3}\vee
c_n}{\sum_{k=1}^nx_{k-l}^2}>\delta\ \hbox{\rm for some}\ n\ge
m\right\} 
$$ 
$$
\subseteq\left\{\frac{|\sum_{k=1}^nx_{k-l}\varepsilon_k|}{(\sum_{k=1}^nx_{k-l}^2)^{2/3}\vee
c_n}>\sqrt{\delta}\ \hbox{\rm for some}\ n\ge m\right\}\cup 
$$ 
$$
\left\{(\sum_{k=1}^nx_{k-l}^2)^{-1/3}\vee
c_n(\sum_{k=1}^nx_{k-l}^2)^{-1}> \sqrt{\delta}\ \hbox{\rm for
some}\ n\ge m\right\}
$$
$$
\subset\left\{|\sum_{k=1}^nx_{k-l}\varepsilon_k|
\left((\sum_{k=1}^nx_{k-l}^2)^{2/3}\vee
c_n\right)^{-1}>\sqrt{\delta}\ \hbox{\rm for some}\ n\ge
m\right\}\cup
$$
$$
\cup\left\{(\sum_{k=1}^nx_{k-l}^2)^{-1}>\delta^{3/2}\ \hbox{\rm
for some}\ n\ge m\right\}\cup
$$
$$
\cup\left\{c_n(\sum_{k=1}^nx_{k-l}^2)^{-1}>\sqrt{\delta}\
\hbox{\rm for some}\ n\ge m\right\} 
$$ 
$$
\subset\left\{|\sum_{k=1}^nx_{k-l}\varepsilon_k|>\delta (\sum_{k=1}^nx_{k-l}^2)^{2/3}\vee
c_n \ \hbox{\rm for some}\ n\ge m\right\}
$$
$$
\cup\left\{c_n(\sum_{k=1}^nx_{k-l}^2)^{-1}>\sqrt{\delta}\wedge\delta^{3/2}\
\hbox{\rm for some}\ n\ge m\right\}.
$$
From here, it follows that
$$
\P_{\theta}
\left(|\sum_{k=1}^nx_{k-l}\varepsilon_k|>\delta\sum_{k=1}^nx_{k-l}^2\
\hbox{\rm for some}\ n\ge m\right)
$$
$$
\le \P_{\theta}
\left(|\sum_{k=1}^nx_{k-l}\varepsilon_k|>\delta\left(\sum_{k=1}^nx_{k-l}^2\right)^{2/3}\vee c_n\
\hbox{\rm for some}\ n\ge m\right)
$$
\begin{equation}\label{3.7}
 +\P_{\theta}
\left\{n^{3/4}(\sum_{k=1}^nx_{k-l}^2)^{-1}>\sqrt{\delta}\wedge\delta^{3/2}\
\hbox{\rm for some}\ n\ge m\right\}.
\end{equation}
$$
$$
By making use of
(\ref{1.1}) and the elementary inequalities, one obtains
$$
\sum_{k=1}^n\varepsilon_k^2=\sum_{k=1}^n(x_{k}-\theta_1x_{k-1}-\ldots-\theta_px_{k-p})^2
$$
$$ \le
(p+1)\left(\sum_{k=1}^nx_{k}^2+\sum_{j=1}^p\theta_j^2\sum_{k=1}^nx_{k-j}^2\right)
$$
$$ \le (p+1)\max_{1\le j\le p}\theta_j^2(\sum_{k=1}^nx_k^2+\sum_{j=1}^p\sum_{k=1}^nx_{k-j}^2)
$$
$$
\le(p+1)^2\mu_k\sum_{k=1}^nx_{k-l}^2
\left(1+\sum_{k=n-l+1}^nx_k^2/\sum_{k=1}^nx_{k-l}^2\right)\,,
$$
where $\mu_k=\sup_{\theta\in K}\|\theta\|^2$.

Therefore
the second summand in the right-hand side of (\ref{3.6}) can be
estimated as
$$
\P_{\theta}\left\{n^{3/4}(\sum_{k=1}^nx_{k-l}^2)^{-1}>\sqrt{\delta}\wedge\delta^{3/2}\
\hbox{\rm for some}\ n\ge m\right\}
$$
$$
\le\P_{\theta}\left\{2(p+1)^2\mu_kn^{3/4}(\sum_{k=1}^n\varepsilon_k^2)^{-1}>\sqrt{\delta}\wedge\delta^{3/2}\
\hbox{\rm for some}\ n\ge m\right\}
$$
$$
+\P_{\theta}\left\{\sum_{k=n-l+1}^nx_k^2/\sum_{k=1}^nx_{k-l}^2\ge 1\ \hbox{\rm for
some}\ n\ge m\right\}.
$$
Combining this and (\ref{3.7}) and
applying Lemma 3.1, and Lemma 2.2 by Lai and Siegmund (see Section 4) yield  (\ref{3.6}).

This completes the proof of Lemma 3.2. $\endproof$

\begin{lemma}\label{L.3.3} Let parameters $\theta_1,\ldots,\theta_p$ in the equation (\ref{1.1})
satisfy Conditions 1-3 and $p\times p$ matrix $M_n$ be given by (\ref{1.3}).
Then, for any compact set $K\subset{\tilde{\Lambda}}_p$ and each $\delta>0$,
$$
\lim_{m\to\infty}\sup_{\theta\in
K}\P_{\theta}\left(\|\frac{M_n}{\sum_{k=1}^nx_{k-1}^2}-
L(\theta_1,\ldots,\theta_p)\|\ge\delta\ \hbox{\rm for some}\ n\ge
m\right)=0,
$$
where $L(\theta_1,\ldots,\theta_p)$ is defined in  (\ref{2.4}).
\end{lemma}
{\sl Proof.} Each diagonal element of the matrix $M_n$ can be expressed through
$\sum_{l=1}^nx_{l-1}^2$ as
\begin{equation}\label{3.8}
<M_n>_{ii}=\sum_{k=i}^nx_{k-i}^2=\sum_{l=1}^{n-i+1}x_{l-1}^2
=\sum_{l=1}^nx_{l-1}^2-\sum_{l=n-i+2}^nx_{l-1}^2,\ 2\le i\le p\,.
\end{equation}
Further it will be observed that each element $<M_n>_{ij},\ 2\le i<j\le p$, of $M_n$
standing above the principal diadonal and below the first row can be expressed through some
element of the first row as
\begin{equation}\label{3.9}
<M_n>_{ij}=\sum_{k=1}^nx_{k-i}x_{k-j}=\sum_{l=1}^nx_{l-1}x_{l-1+i-j}
-\sum_{l=n-i+2}^{n}x_{l-1}x_{l-1+i-j}
\end{equation}
$$
=<M_n>_{i,j-i+1}-\sum_{t=n-i+2}^nx_{t-1+i-j} \,.
$$
Now we derive the equations relating the elements of the first row, that is, $<M_n>_{1s},\ 2\le s\le p$.
Making use of the equation (\ref{1.1}), one gets
$$
<M_n>_{1s}=\sum_{k=1}^nx_{k-1}x_{k-s}
=\sum_{k=1}^n\left(\sum_{l=1}^p\theta_lx_{k-l-1}+\varepsilon_{k-1}\right)x_{k-s}
$$
$$
=\sum_{l=1}^p\theta_l\sum_{k=1}^nx_{k-l-1}x_{k-s}+\sum_{k=1}^nx_{k-s}\varepsilon_{k-1}\,.
$$ 
Since
$$
\sum_{k=1}^nx_{k-l-1}x_{k-s}=
\left\{
\begin{array}{cc}
<M_n>_{ss}, & \hbox{\rm if}\ l=s-1,\\
<M_n>_{s,l+1}, & \hbox{\rm if}\ l\ge s,\\
<M_n>_{l+1,s},  & \hbox{\rm if}\ l\le s-2,
\end{array}
\right.
$$
this implies the following system of equations
$$
<M_n>_{12}=\theta_1<M_n>_{22}+\sum_{l=2}^p\theta_l<M_n>_{2,l+1}+\sum_{k=1}^nx_{k-2}\varepsilon_{k-1}\,,
$$
$$
<M_n>_{1s}=\sum_{l=2}^{s-2}\theta_l<M_n>_{l+1,s}+\theta_{s-1}<M_n>_{ss}
$$
$$
+\sum_{l=s}^p\theta_l<M_n>_{s,l+1}+\sum_{k=1}^nx_{k-s}\varepsilon_{k-1},\ 3\le s\le p.
$$
Taking into account (\ref{3.8}),(\ref{3.9}) one can represent this system as
\begin{equation}\label{3.10}
<M_n>_{12}=\theta_1\sum_{k=1}^nx_{k-1}^2+\sum_{l=2}^p\theta_l<M_n>_{1,l}+\eta_{1,2}(n),
\end{equation}
$$
<M_n>_{1s}=\sum_{l=1}^{s-2}\theta_l<M_n>_{1,s-l}+\theta_{s-1}\sum_{k=1}^nx_{k-1}^2
+\sum_{l=s}^p\theta_l<M_n>_{1,l-s+2}+\eta_{1,s}(n),
$$
where $ 3\le s\le p,$
$$
\eta_{1,2}(n)=-\theta_1x_{n-1}^2-\sum_{l=2}^p\theta_lx_{n-1}x_{n-l}+\sum_{k=1}^nx_{k-2}\varepsilon_{k-1},
$$
$$
\eta_{1,s}(n)=-\sum_{l=1}^{s-2}\theta_l\sum_{t=n-l+1}^nx_{t-1}x_{t+l-s}-\theta_{s-1}
\sum_{l=n-s+2}^nx_{l-1}^2
$$
$$
-\sum_{l=s}^p\theta_l\sum_{t=n-s+2}^nx_{t-1}x_{t+s-l-2}+\sum_{k=1}^nx_{k-s}\varepsilon_{k-1}.
$$
Denote
$$
z_i(n)=\frac{<M_n>_{1,i+1}}{\sum_{k=1}^nx_{k-1}^2},\
{\tilde{\eta}}_{1,i+1}=\frac{\eta_{1,i+1}(n)}{\sum_{k=1}^nx_{k-1}^2},\ i=\overline{1,p-1}.
$$ 
Then the system of equations (\ref{3.10}) takes the form
\begin{equation}\label{3.11}
\left\{
\begin{array}{l}
z_1(n)-\sum_{l=2}^p\theta_lz_{l-1}(n)=\theta_1+{\tilde{\eta}}_{1,2},\\
-\sum_{k=1}^{j-1}\theta_{j-k}z_k(n)+z_j(n)-\sum_{k=1}^{p-j}\theta_{k+j}z_k(n)=\theta_j+{\tilde{\eta}}_{1,j+1},\\
j=\overline{2,p-1}
\end{array}
\right.
\end{equation}
In virtue of Lemmas 3.1,3.2, for any compact set $K\subset{\tilde{\Lambda}}_p$ and $\delta>0$, one has
$$
\lim_{m\to\infty}\sup_{\theta\in K}\P_{\theta}(\max_{2\le s\le p}|{\tilde{\eta}}_{1,s}(n)|>\delta\
 \hbox{\rm for some}\ n\ge m)=0.
$$
From here and the Condition 3, which holds for each vector $\theta\in K$, it follows that the solution
 of the system (\ref{3.11}) converges, as $n\to\infty$, to the unique solution of system (\ref{2.3})
 uniformly in $\theta\in K$, that is, 
 $$
 \lim_{m\to\infty}\sup_{\theta\in K}\P_{\theta}(\max_{1\le i\le p-1}|z_i(n)-\kappa_i(\theta)|>\delta\
 \hbox{\rm for some}\ n\ge m)=0.
 $$
 This, in view of (\ref{3.9}) and Lemma 3.1, implies the desired convergence of the remaining elements
 of the matrix $M_n$. Hence Lemma 3.3. $\endproof$

\begin{lemma}\label{L.3.4}
Let $M_n,\tau(h)$ and $L(\theta)=L(\theta_1,\ldots,\theta_p)$ be given by (\ref{1.4}),
(\ref{2.2}) and (\ref{2.4}),  respectively.  Then, for any compact set $K \subset
{\stackrel{\circ} {\Lambda}}_p$ and $\delta>0$,
\begin{equation}\label{3.12}
\lim_{h\to\infty}\sup_{\theta\in
K}\P_{\theta}\left(\|\frac{M_{\tau(h)}}{h}-L(\theta_1,\ldots,\theta_p)\|>\delta\right)=0,
\end{equation}
where ${\stackrel{\circ} {\Lambda}}_p$ is defined in Condition 3.
\end{lemma}
{\sl Proof.} By making use of the equality
$$
\frac{M_{\tau(h)}}{h}-L(\theta)=\frac{M_{\tau(h)}}{\sum_{k=1}^{\tau(h)}x_{k-1}^2}
-L(\theta)+M_{\tau(h)}\left(\frac{1}{h}-\frac{1}{\sum_{k=1}^{\tau(h)}x_{k-1}^2}\right)
$$
one gets the estimate
$$ \|\frac{M_{\tau(h)}}{h}-L(\theta)\|
\le\|\frac{M_{\tau(h)}}{\sum_{k=1}^{\tau(h)}x_{k-1}^2}-L(\theta)\|
$$
$$
+\left(\|M_{\tau(h)}\|/\sum_{k=1}^{\tau(h)}x_{k-1}^2\right)
\left(\sum_{k=1}^{\tau(h)}x_{k-1}^2-h\right)/h.
$$
Therefore
\begin{equation}\label{3.13}
\left\{\|\frac{M_{\tau(h)}}{h}-L(\theta)\|>\delta\right\}
\subset\left\{\|\frac{M_{\tau(h)}}{\sum_{k=1}^{\tau(h)}x_{k-1}^2}-L(\theta)\|
>\delta/2\right\}
\end{equation}
$$
\cup\left\{\|\frac{M_{\tau(h)}}{\sum_{k=1}^{\tau(h)}x_{k-1}^2}\|\frac{\left(\sum_{k=1}^{\tau(h)}x_{k-1}^2-h\right)}{h}>\delta/2\right\}.
$$
Further one has the inclusions
$$
\left\{\|\frac{M_{\tau(h)}}{\sum_{k=1}^{\tau(h)}x_{k-1}^2}-L(\theta)\|>\delta/2\right\}
\subseteq\{\tau(h)\le m\}
$$
$$
\cup\left\{\|\frac{M_{n}}{\sum_{k=1}^{n}x_{k-1}^2}-L(\theta)\|>\delta/2\
\hbox{\rm for some}\ n\ge m\right\},
$$
$$
\left\{\|\frac{M_{\tau(h)}}{\sum_{k=1}^{\tau(h)}x_{k-1}^2}\|\left(\sum_{k=1}^{\tau(h)}x_{k-1}^2-h\right)/h>\delta/2\right\}
\subset\{\tau(h)\le m\}
$$
$$
\cup\left\{\|\frac{M_{n}}{\sum_{k=1}^{n}x_{k-1}^2}-L(\theta)\|>\delta/2\
\hbox{\rm for some}\ n\ge m\right\},
$$
$$
\left\{\|\frac{M_{\tau(h)}}{\sum_{k=1}^{\tau(h)}x_{k-1}^2}\|\left(\sum_{k=1}^{\tau(h)}x_{k-1}^2-h\right)/h>\delta/2\right\}
\subset\{\tau(h)\le m\}
$$
$$
\cup\left\{\frac{\|M_n\|}{\sum_{k=1}^{n}x_{k-1}^2}\cdot\frac{x_{n-1}^2}{\sum_{k=1}^{n-1}x_{k-1}^2}>\delta/2\
\hbox{\rm for some}\ n\ge m\right\}.
$$
From here and (\ref{3.13}), it follows that
\begin{equation}\label{3.14}
\P_{\theta}\left\{\|\frac{M_{\tau(h)}}{h}-L(\theta)\|>\delta\right\}
\le 2\P_{\theta}\{\tau(h)\le m\}
\end{equation}
$$
+\P_{\theta}\left\{\|\frac{M_{n}}{\sum_{k=1}^{n}x_{k-1}^2}-L(\theta)\|>\delta/2\
\hbox{\rm for some}\ n\ge m\right\}
$$
$$
+\P_{\theta}\left\{\frac{\|M_n\|}{\sum_{k=1}^{n}x_{k-1}^2}\cdot\frac{x_{n-1}^2}{\sum_{k=1}^{n-1}x_{k-1}^2}>\delta/2\
\hbox{\rm for some}\ n\ge m\right\}. $$
By the definition of $\tau(h)$
in (\ref{2.2})
$$
\{\tau(h)<m\}=\left\{\sum_{k=1}^{m}(x_{k-1}^2+\cdots+x_{k-p}^2)>h\right\}
$$
$$
=\left\{\sum_{k=1}^{m}(x_{k-1}^2+\cdots+x_{k-p}^2)>h,\ \max_{1\le
j\le m}(x_{j-1}^2+\cdots+x_{j-p}^2)<l\right\}
$$
$$
+\left\{\sum_{k=1}^{m}(x_{k-1}^2+\cdots+x_{k-p}^2)>h,\ \max_{1\le j\le
m}(x_{j-1}^2+\cdots+x_{j-p}^2)\ge l\right\}
$$
\begin{equation}\label{3.15}
\subset\{ml>h\}\cup\cup_{j=1}^m\{(x_{j-1}^2+\cdots+x_{j-p}^2)\ge l\}.
\end{equation}
This yields
$$
\P_{\theta}\{\tau(h)<m\}\le
I_{(ml>h)}+\sum_{k=1}^{m}\P_{\theta}\{x_{k-1}^2+\cdots+x_{k-p}^2\ge l\}.
$$
Consider the last term in (\ref{3.14}).
By the inequality
$$ \frac{\|M_n\|}{\sum_{k=1}^{n}x_{k-1}^2}\le
\|\frac{M_{n}}{\sum_{k=1}^{n}x_{k-1}^2}-L(\theta)\|+\|L(\theta)\|
$$
one has
$$
\P_{\theta}\left\{\frac{\|M_n\|}{\sum_{k=1}^{n}x_{k-1}^2}
\cdot\frac{x_{n-1}^2}{\sum_{k=1}^{n-1}x_{k-1}^2}>\delta/2\
\hbox{\rm for some}\ n\ge m\right\}
$$
\begin{equation}\label{3.16}
\le\P_{\theta}\left\{\|\frac{M_{n}}{\sum_{k=1}^{n}x_{k-1}^2}-L(\theta)\|>\sqrt{\delta/4}\
 \hbox{\rm for some}\ n\ge m\right\}
\end{equation}
$$
+\P_{\theta}\left\{L_k^*\frac{x_{n-1}^2}{\sum_{k=1}^{n-1}x_{k-1}^2}>\sqrt{\delta/4}\
\hbox{\rm for some}\ n\ge m\right\}
$$
$$
+\P_{\theta}\left\{L_k^*\frac{x_{n-1}^2}{\sum_{k=1}^{n-1}x_{k-1}^2}>\delta/4\\\
\hbox{\rm for some}\ n\ge m\right\},
$$
where
$L_k^*=\sup_{\theta\in K}\|L(\theta)\|$.

Combining (\ref{3.14})-(\ref{3.16}) yields
$$
\sup_{\theta\in
K}\P_{\theta}\left(\|\frac{M_{\tau(h)}}{h}-L(\theta)\|>\delta\right)
$$
$$
\le 2I_{(ml>h)}+2\sum_{k=1}^{m}\sup_{\theta\in
K}\P_{\theta}\{\|X_{k-1}\|^2\ge l\}
$$
$$
+\P_{\theta}\left\{\|\frac{M_{n}}{\sum_{k=1}^{n}x_{k-1}^2}-L(\theta)\|>
\frac{1}{2}(\delta\wedge\sqrt{\delta})\\\ \hbox{\rm for some}\ n\ge
m\right\}
$$
$$
+2\P_{\theta}\left\{L_k^*\frac{x_{n-1}^2}{\sum_{k=1}^{n-1}x_{k-1}^2}>
\frac{1}{2}(\delta\wedge\sqrt{\delta})\\\ \hbox{\rm for some}\ n\ge
m\right\}.
$$
 Limiting $h\to\infty,\ l\to\infty,\ m\to\infty$ and
taking into account Lemma 3.3, one comes to (\ref{3.12}). Hence Lemma
3.4. $\endproof$

\begin{lemma}\label{L.3.5.} Let $x_k$ and $\tau(h)$ be defined by (\ref{1.1}) and
 (\ref{2.2}). Then for any compact set $K\subset{\stackrel{\circ}{\Lambda}}_0$ and $\delta>0$
\begin{equation}
\label{3.17}
    \lim_{h\to\infty}\sup_{\theta\in K}\P_{\theta}\left(x_{\tau(h)-1}^2/
    \sum_{k=1}^{\tau(h)-1}x_{k-1}^2>\delta\right)=0.
\end{equation}
\end{lemma}
 {\sl Proof.} In view of the inclusion
$$
    \left\{\frac{x_{\tau(h)-1}^2}{\sum_{k=1}^{\tau(h)-1}x_{k-1}^2}>\delta\right\}
    \subset
    \left\{\tau(h)\leq m\right\}\bigcup\left\{\frac{x_{n}^2}{\sum_{k=1}^{n-1}x_{k-1}^2}
    >\delta \\\  \hbox{\rm for some} \; n\geq m \right\}.
    $$
one has
    $$
   \sup_{\theta\in K}\P_{\theta}\left\{\frac{x_{\tau(h)-1}^2}{\sum_{k=1}^{\tau(h)-1}
   x_{k-1}^2}>\delta\right\}\leq I_{\left(ml>h\right)}+\sum_{j=1}^{m}\sup_{\theta
   \in K}\P_{\theta}\left\{\|X_{j-1}\|^2 \geq l\right\}
   $$
   $$
   +\sup_{\theta\in K}
   \P_{\theta}\left\{\frac{x_n^2}{\sum_{k=1}^{n-1}x_{k-1}^2}>\delta \\\ \hbox{\rm for some} \; n\geq m\right\}.
$$
 Limiting $ h\to\infty$, $l\to\infty$, $m\to\infty$ and applying
Lemma 3.1 lead to (\ref{3.17}). This completes the proof of Lemma
3.5. $\endproof$

\sect{Appendix.}

In this Section we cite the probabilistic result from the paper of Lai and
Siegmund (1983) used in Section 3 and give the proof of Theorem 2.2.

{\bf Lemma 2.2} (by Lai and Siegmund (1983)).
Let $({\cF}_n)_{n\ge 0}$ be a filtration on a measurable space $(\Omega,\cF),\ (x_n)_{n\ge 0}$
and $(\varepsilon_n)_{n\ge 0}$ be sequences of random variables adapted to $({\cF}_n)_{n\ge 0}$.
Let $(\P_{\theta}, \theta\in\Theta)$ be a family of probability measures on $(\Omega,\cF)$ such
that under every $\P_{\theta}\  \varepsilon_n$ is independent of ${\cF}_{n-1}$ for each $n\ge 1$. Then, for each $\gamma>1/2, \delta>0$, and increasing sequence of positive constants $c_n\to\infty$,
$$
\sup_{\theta\in\Theta}\P_{\theta}\left(|\sum_{i=1}^nx_{i-1}\varepsilon_i|\ge\delta
\max(c_n,(\sum_{i=1}^nx_{i-1}^2)^{\gamma})\ \hbox{\rm for some}\ n\ge m\right)\to 0\,,
$$
as $m\to\infty$.

{\sl  Proof of Theorem 2.2}. Assertion (\ref{2.20}) easely follows from Lemma
3.12 in [6]. For the points $\theta$ belonging to $\partial\Lambda_p$ we decompose the original time
series (\ref{1.2}) into several components depending on the number of the roots of the
characteristic polynomial (\ref{1.5}) lying on the unit circle and their values. To this end,
the characteristic polynomial (\ref{1.5}) is represented as
$$
\cP(z)=(z+1)^{\delta_1}(z-1)^{\delta_2}(z^2-2z\cos\phi\,+\,1)^{\delta_3}\varphi(z),
$$
where $\delta_i$ are either zero or $1$ with $\delta_1+\delta_2+\delta_3\ge 1,\ \varphi(z)$
is the polynomial of order $r=p-\delta_1-\delta_2-2\delta_3$ which has all roots inside the
 unit circle. Assuming (without loss of generality) that $r\ge 1$, one has
 $$
 \varphi(z)= z^r+\beta_1z^{r-1}+\ldots+\beta_r\,.
 $$
 By applying the backshift operator $q^{-1}$ (i.e. $q^{-1}x_n=x_{n-1}$) one can write down
 (\ref{1.1}) as
\begin{equation}\label{4.1}
q^{-\delta_1}(q+1)^{\delta_1}q^{-\delta_2}(q-1)^{\delta_2}q^{-2\delta_3}(q^2-2q\cos\phi+1)^{\delta_3}
q^{-r}\varphi(q)x_n=\varepsilon_n
\end{equation}
Let $\theta\in\Gamma_1(\rho)$. Then this equation, in view of (\ref{2.19}), takes the form
$$
q^{-1}(q+1)q^{-p+1}\varphi(q)x_n=\varepsilon_n\,.
$$
Denote 
$$
u_n=q^{-p+1}\varphi(q)x_n\,,
$$
$$
v_n=q^{-1}(q+1)x_n\,,
$$
that is
\begin{equation}\label{4.2}
\begin{array}{l}
u_n=x_n+\beta_1x_{n-1}+\ldots+\beta_{p-1}x_{n-p+1}\,,\\
v_n=x_n+x_{n-1}\,.
\end{array}
\end{equation}
 Introducing the vector $V_n=(v_n,\ldots,v_{n-p+1})^{\prime}$
and the matrix
\begin{equation}\label{4.3}
Q=
\left(
\begin{array}{ccccc}
1 & \beta_1 & \beta_2 & \ldots & \beta_{p-1} \\
1 & 1       & 0       & \ldots & 0           \\
0 & 1       & 1       & \ldots & 0           \\
\vdots      & {}      & {}                   \\
0 &         & \ldots  & 1      & 1           \\
\end{array}
\right)
\end{equation}
one can rewrite equations (\ref{4.2}) in the vector form
\begin{equation}\label{4.4}
\left(
\begin{array}{c}
U_n\\
V_n
\end{array}
\right)
=Q_1X_n,\\\ X_n=(x_n,\ldots, x_{n-p+1})^{\prime}\,.
\end{equation}
The processes $u_n$ and $v_n$ satisfy the equations
$$
u_n=-u_{n-1}+\varepsilon_n,\\\ v_n+\beta_1v_{n-1}+\ldots+\beta_{p-1}v_{n-p+1}=\varepsilon_n\,.
$$
Substituting (\ref{4.4}) in (\ref{1.4}) yields
\begin{equation}\label{4.5}
\frac{tr\,M_n}{n^2}=tr(Q_1Q_1^{\prime})^{-1}
\left(
\begin{array}{cc}
n^{-2}\sum_{k=1}^nu_{k-1}^2 & n^{-2}\sum_{k=1}^nu_{k-1}V_{k-1}^{\prime}\\
n^{-2}\sum_{k=1}^nu_{k-1}V_{k-1} & n^{-2}\sum_{k=1}^nV_{k-1}V_{k-1}^{\prime}
\end{array}
\right).
\end{equation}
By Theorem 3.4.2. in Chan and Wei (1988) 
$$
n^{-2}\sum_{k=1}^nu_{k-1}V_{k-1}^{\prime}\stackrel{\P}{\to}0,\ \hbox{\rm as}\ n\to \infty.
$$
Since the process $V_n$ is stable
$$
\lim_{n\to\infty}n^{-2}\sum_{k=1}^nV_{k-1}V_{k-1}^{\prime}=0\\\ \hbox{\rm a.s.}
$$
By Donsker's theorem
$$
n^{-2}\sum_{k=1}^{[nt]}u_{k-1}^2\stackrel{\mathcal{L}}{\Longrightarrow}\sigma^2\int_0^t\,W_1^2(s)ds,\ 0\le t\le 1\,,
$$
as $n\to \infty$. By making use of these limiting relations in (\ref{4.5}), one get
\begin{equation}\label{4.6}
\frac{tr\,M_{[nt]}}{n^2}\stackrel{\mathcal{L}}{\Longrightarrow}\kappa_{11}\int_0^t\,W_1^2(s)ds,\ 0\le t\le 1\,,
\end{equation}
where
$$
\kappa_{11}=\sigma^2<(Q_1Q_1^{\prime})^{-1}>_{11}\,.
$$
Now by definition of $\tau(h)$ in (\ref{2.2}) one has 
$$
\P_{\theta}(\frac{\tau(h)}{b_1\sqrt{h}}\le t)=\P_{\theta}(tr\,M_{[tb_1\sqrt{h}]}\ge h)
=\P_{\theta}(\frac{tr\,M_{[tb_1\sqrt{h}]}}{b_1^2h}b_1^2\ge 1)\,.
$$
This and (\ref{4.6}) imply the validity of (\ref{2.21}) for $\theta\in\Gamma_1(\rho)$ with
\begin{equation}\label{4.7}
b_1^2=1/\kappa_{11}\,.
\end{equation}
By a similar argument, one check (\ref{2.21}) for $\theta\in\Gamma_2(\rho)$ with
\begin{equation}\label{4.8}
b_2^2=1/(\sigma^2<(Q_2Q_2^{\prime})^{-1}>_{11})\,,
\end{equation} 
where
$$
Q_2=
\left(
\begin{array}{ccccc}
1 & \beta_1 & \beta_2 & \ldots & \beta_{p-1} \\
1 & -1      &  0      & \ldots & 0   \\
0 & 1      &  -1       & \ldots & 0  \\
\vdots & {}     &  {}      &  {}      &0 \\
0 &   {}     & \ldots   &  1       & -1
\end{array}
\right).
$$
Assume that $\theta\in\Gamma_3(\rho)$. Then using the equation
$$
q^{-2}(q^2-2q\cos\phi-1)q^{-p+2}\varphi(q)x_n=\varepsilon_n\,,
$$
we decompose the process (\ref{1.2}) into two processes
$$
U_n=(u_n,u_{n-1})^{\prime},\\\ V_n=(v_n,\ldots,v_{n-p+3})^{\prime}\,,
$$
which obey the equations
\begin{equation}\label{4.9}
U_n=
\left(
\begin{array}{cc}
2\cos\phi & -1  \\
1 & 0  
\end{array}
\right)
U_{n-1}+
\left(
\begin{array}{c}
\varepsilon_n  \\
 0  
\end{array}
\right)\,,
\end{equation}
$$
v_n+\beta_1v_{n-1}+\ldots+\beta_{p-2}v_{n-p+2}=\varepsilon_n\,.
$$
These processes are related with $X_n=(x_n,\ldots,x_{n-p+1})^{\prime}$ by the following
transformation
$$
\left(
\begin{array}{c}
U_n  \\
 V_n  
\end{array}
\right)
=Q_3X_n
$$
with
$$
Q_3=
\left(
\begin{array}{cccccc}
1 & \beta_1    & \ldots &  {}    &\beta_{p-2}& 0  \\
0 & 1          & \beta_1& \ldots & {}        & \beta_{p-2}\\
1 & -2\cos\phi & 1      & 0      & \ldots    & 0          \\
0 & 1          &-2\cos\phi& 1    & 0         &\ldots      \\
\vdots& {}     & {}       & {}   & {}        & {}         \\
0 & \ldots     &  0       & 1    & -2\cos\phi& 1
\end{array}
\right).
$$
Further, by the same argument, one shows that
$$
n^{-2}\,tr\,M_{[nt]}\stackrel{\mathcal{L}}{\Longrightarrow}tr(Q_3Q_3^{\prime})^{-1}
\left(
\begin{array}{cc}
H_t & 0  \\
 0  & 0 
\end{array}
\right)\,,
0\le t\le 1,
$$
where
$$
H_t=\frac{\sigma^2}{4\sin^2\phi}\int_0^t(W_1^2(s)+W_2^2(s))ds
\left(
\begin{array}{cc}
1        & \cos\phi  \\
\cos\phi & 1 
\end{array}
\right)\,,
$$
that is,
$$
n^{-2}\,tr\,M_{[nt]}\stackrel{\mathcal{L}}{\Longrightarrow}\frac{r^2}{4\sin^2\phi}\int_0^t(W_1^2(s)+W_2^2(s))ds\,,
$$
where
$$
r^2=\sigma^2\,tr(Q_3Q_3^{\prime})^{-1}
\left(
\begin{array}{cc}
1        & \cos\phi  \\
\cos\phi & 1 
\end{array}
\right)\,.
$$
From here one comes to (\ref{2.21}) with
\begin{equation}\label{4.10}
b_3=(2\sin\phi)/r\,.
\end{equation}
Now assume that $\theta\in\Gamma_4(\rho)$. In this case the process (\ref{1.2}) is decomposed
into two vector processes 
$$
U_n=(u_n,u_{n-1})^{\prime},\ V_n=(v_n,\ldots,v_{n-p+3})^{\prime}\,,
$$ 
which are defined by the formulae
$$
u_n=x_n+\beta_1x_{n-1}+\ldots+\beta_{p-2}x_{n-p+2},\ v_n=x_n-x_{n-2}
$$
and satisfy the equations
$$
u_n=u_{n-2}+\varepsilon_n,\
v_n+\beta_1v_{n-1}+\ldots+\beta_{p-2}v_{n-p+2}=\varepsilon_n\,.
$$
These processes are related to the original process (\ref{1.2}) by the 
transformation
\begin{equation}\label{4.11}
\left(
\begin{array}{c}
U_n  \\
 V_n  
\end{array}
\right)
=Q_4X_n\,,
\end{equation}
where
$$
Q_4=
\left(
\begin{array}{cccccc}
1 & \beta_1    & \ldots &  {}    &\beta_{p-2}& 0  \\
0 & 1          & \beta_1& \ldots & {}        & \beta_{p-2}\\
1 & 0          & -1      & 0      & \ldots    & 0          \\
\vdots& {}     & {}       & {}   & {}        & {}         \\
0 & \ldots     &  0       & 1    & 0         & -1
\end{array}
\right)\,.
$$
Further we represent $U_n$ as
\begin{equation}\label{4.12}
U_n=T^{-1}
\left(
\begin{array}{c}
y_n  \\
 z_n  
\end{array}
\right)\,,
\end{equation}
where processes $y_n$ and $z_n$ satisfy the equations
\begin{equation}\label{4.13}
y_n=y_{n-1}+\varepsilon,\ z_n=-z_{n-1}+\varepsilon\,,
\end{equation}
$$
T=
\left(
\begin{array}{cc}
1 & 1  \\
1 & -1  
\end{array}
\right)\,.
$$
By making use of (\ref{4.11})-(\ref{4.12}) one gets
\begin{equation}\label{4.14}
\frac{tr\,M_n}{n^2}=tr(Q_4Q_4^{\prime})^{-1}
\left(
\begin{array}{cc}
n^{-2}\sum_{k=1}^nU_{k-1}U_{k-1}^{\prime} & n^{-2}\sum_{k=1}^nU_{k-1}V_{k-1}^{\prime}  \\
n^{-2}\sum_{k=1}^nV_{k-1}U_{k-1}^{\prime} & n^{-2}\sum_{k=1}^nV_{k-1}V_{k-1}^{\prime}  
\end{array}
\right)\,,
\end{equation}
$$
n^{-2}\sum_{k=1}^nU_{k-1}U_{k-1}^{\prime}=T^{-1}
\left(
\begin{array}{cc}
n^{-2}\sum_{k=1}^ny_{k-1}^2 & n^{-2}\sum_{k=1}^ny_{k-1}z_{k-1}  \\
n^{-2}\sum_{k=1}^ny_{k-1}z_{k-1} & n^{-2}\sum_{k=1}^nz_{k-1}^2  
\end{array}
\right)(T^{-1})^{\prime}\,.
$$
By Theorems 3.4.1 and 2.3 in Chan and Wei (1988), one has, as $n\to\infty$,
$$
n^{-2}\sum_{k=1}^nU_{k-1}V_{k-1}^{\prime}\stackrel{\P}{\to}0,\\\
n^{-2}\sum_{k=1}^ny_{k-1}z_{k-1}\stackrel{\P}{\to}0\,,
$$
$$
\sigma^{-2}\left(n^{-2}\sum_{k=1}^{[nt]}y_{k-1}^2, n^{-2}\sum_{k=1}^{[nt]}z_{k-1}^2\right)
\stackrel{\mathcal{L}}{\Longrightarrow}\left(J_1(W_1;t), J_1(W_2;t)\right),
$$
where $0\le t\le 1$. From here and (\ref{4.11}), it follows that
\begin{equation}\label{4.15}
n^{-2}tr\,M_{[nt]}\stackrel{\mathcal{L}}{\Longrightarrow}\sigma^2tr(Q_4Q_4^{\prime})^{-1}{\cD}_t\,,
\end{equation}
where
$$
{\cD}_t=
\left(
\begin{array}{cc}
S_t & 0  \\
0 & 0  
\end{array}
\right)\,,
S_t=T^{-1}
\left(
\begin{array}{cc}
\int_0^tW_1^2(s)ds & 0  \\
0 & \int_0^tW_2^2(s)ds  
\end{array}
\right)
(T^{-1})^{\prime}\,.
$$
It easy to check that
\begin{equation}\label{4.16}
\sigma^2tr(Q_4Q_4^{\prime})^{-1}{\cD}_t=\int_0^t\left(r_1W_1^2(s)+r_2W_2^2(s)\right)ds\,,
\end{equation}
where
$$
r_1=\frac{\sigma^2}{4}\sum_{i=1}^2\sum_{j=1}^2<(Q_4Q_4^{\prime})^{-1}>_{ij},
r_2=\frac{\sigma^2}{4}\sum_{i=1}^2\sum_{j=1}^2(-1)^{i+j}<(Q_4Q_4^{\prime})^{-1}>_{ij}\,.
$$
Denote
\begin{equation}\label{4.17}
b_4^2=1/\mu_1,\\\ \mu_1=r_2/r_1\,.
\end{equation}
It remains to note that (\ref{4.15}),(\ref{4.16}) imply
$$
\P_{\theta}(\tau(h)\le tb_4\sqrt{h})\to\P_{\theta}(J_3(W_1,W_2)\ge 1),\ \hbox{\rm as}\ h\to\infty\,,
$$
where $J_3$ is defined in (\ref{2.18}).

Assume that $\theta\in\Gamma_5(\rho)$. Then equation (\ref{4.1}) has the form
$$
q^{-1}(q+1)q^{-2}(q^2-2q\cos\phi+1)q^r\varphi(q)x_n=\varepsilon_n,\ r=p-3.
$$
Decompose $x_n$ into three processes
\begin{equation}\label{4.18}
\begin{array}{l}
u_n=q^{-2}(q^2-2q\cos\phi+1)q^{-r}\varphi(q)x_n
=x_n+\gamma_1x_{n-1}+\ldots+\gamma_{p-1}x_{n-p+1}\,,\\
v_n=q^{-1}(q+1)q^{-r}\varphi(q)x_n=x_n+f_1x_{n-1}+\ldots+f_{p-2}x_{n-p+2}\,,\\
 w_n=q^{-1}(q+1)q^{-2}(q^2-2q\cos\phi+1)x_n=x_n+t_1x_{n-1}+t_2x_{n-2}+t_3x_{n-3}\,,
 \end{array}
\end{equation}
where $\gamma_i, f_j, t_k$ are the coefficients of the corresponding polynomials. These
processes satisfy the following equations
$$
\begin{array}{l}
u_n=-u_{n-1}+\varepsilon_n,\\
v_n=2v_{n-1}\cos\phi-v_{n-2}+\varepsilon_n,\\
w_n=-\beta_1 w_{n-1}-\ldots-\beta_{p-3} w_{n-p+3}+\varepsilon_n\,.
\end{array}
$$
Introducing vectors $V_n=(v_n,v_{n-1})^{\prime},\ W_n=(w_n,w_{n-1},\ldots,w_{n-p+4})^{\prime}$
and the matrix 
\begin{equation}\label{4.19}
Q_5=
\left(
\begin{array}{cccccccc}
1 & \gamma_1 & \ldots & {} & {} & {} & {}& \gamma_{p-1}\\
1 & f_1      & \ldots & {} & {} & {} & f_{p-2} & 0 \\
0 & 1 & f_1 & \ldots & {} & {} & {} & f_{p-2}\\
1 & t_1 & t_2 & t_3 & 0 & \ldots & {} & 0 \\
0 & 1 & t_1 & t_2 & t_3 & 0 & \ldots & 0 \\
\ldots & {} &{}&{}&{}&{}&{}&{}\\
0 & \ldots &{}&0 & 1 & t_1 & t_2 & t_3
\end{array}
\right)\,,
\end{equation}
we write down (\ref{4.18}) as
\begin{equation}\label{4.20}
\left(
\begin{array}{c}
U_n\\
V_n\\
W_n
\end{array}
\right)
=Q_5X_n\,.
\end{equation}
From here and (\ref{1.4}) 
\begin{equation}\label{4.21}
n^{-2}tr\,M_{[nt]}=tr(Q_5Q_5^{\prime})^{-1}C_{[nt]}/n^2\,,\ 0\le t\le 1,
\end{equation}
where
$$
C_n=
\left(
\begin{array}{ccc}
\sum_{k=1}^nu_{k-1}^2 & \sum_{k=1}^nU_{k-1}V_{k-1}^{\prime} & \sum_{k=1}^nU_{k-1}W_{k-1}^{\prime}\\
\sum_{k=1}^nV_{k-1}u_{k-1} & \sum_{k=1}^nV_{k-1}V_{k-1}^{\prime} &\sum_{k=1}^nV_{k-1}W_{k-1}^{\prime}\\
\sum_{k=1}^nW_{k-1}u_{k-1}& \sum_{k=1}^nW_{k-1}V_{k-1}^{\prime}& \sum_{k=1}^nW_{k-1}W_{k-1}^{\prime}
\end{array}
\right).
$$
By Theorems 3.4.1, 3.4.2 in Chan and Wei (1988),
\begin{equation}\label{4.22}
\begin{array}{l}
n^{-2}\sum_{k=1}^nV_{k-1}u_{k-1}\stackrel{\P}{\to}0\,,\\
n^{-2}\sum_{k=1}^nW_{k-1}u_{k-1}\stackrel{\P}{\to}0\,,\\
n^{-2}\sum_{k=1}^nW_{k-1}V_{k-1}^{\prime}\stackrel{\P}{\to}0\,,
\end{array}
\end{equation}
as $n\to\infty$. Due to the stability of the process $W_k$,
\begin{equation}\label{4.23}
\lim_{n\to\infty}n^{-2}\sum_{k=1}^nW_{k-1}W_{k-1}^{\prime}=0\ \hbox{\rm a.s.}
\end{equation}
By the Theorem 2.3 therein
\begin{equation}\label{4.24}
\left(n^{-2}\sum_{k=1}^{[nt]}u_{k-1}^2,\ n^{-2}\sum_{k=1}^{[nt]}V_{k-1}V_{k-1}^{\prime}\right)
\stackrel{\mathcal{L}}{\Longrightarrow}\sigma^2\left(J_2(W_1,W_2;t)H,\ J_1(W_3;t)\right)\,,
\end{equation}
where $J_1$ and $J_2$ are defined in (\ref{2.18}),
$$
H=\frac{1}{4\sin^2\phi}
\left(
\begin{array}{cc}
1 & \cos\phi \\
\cos\phi & 1
\end{array}
\right).
$$
Limiting $n\to\infty$ in (\ref{4.21}) and taking into account (\ref{4.22})-(\ref{4.24}) one gets
\begin{equation}\label{4.25}
n^{-2}tr\,M_{[nt]}\stackrel{\mathcal{L}}{\Longrightarrow}m(t)\ \hbox{\rm as}\ n\to\infty\,,
\end{equation}
where
\begin{equation}\label{4.26}
m(t)=\sigma^2tr(Q_5Q_5^{\prime})^{-1}
\left(
\begin{array}{cc}
S_t & 0 \\
0   & 0
\end{array}
\right),
\end{equation}
$$
S_t=
\left(
\begin{array}{cc}
J_1(W_3;t) & 0 \\
0   & J_2(W_1,W_2;t)H
\end{array}
\right).
$$
By easy calculation in (\ref{4.26}) one finds
$$
m(t)=\kappa_{11}J_1(W_3;t)+\frac{J_2(W_1,W_2;t)}{4\sin^2\phi}(\kappa_{22}+\kappa_{33}+(\kappa_{23}+\kappa_{32})\cos\phi)\,,
$$
where $\kappa_{ij}=\sigma^2<(Q_5Q_5^{\prime})^{-1}>_{ij}$. 

Now by the definition of $\tau(h)$ in (\ref{2.2}) and (\ref{4.25}) we obtain
$$
\P(\tau(h)\le b_5\sqrt{h}t)=\P(tr\,M_{[ tb_5\sqrt{h}]}\ge h)\,,
$$
$$
\P(\frac{tr\,M_{[ tb_5\sqrt{h}]}}{hb_5^2}b_5^2\ge 1)\to\P(m(t)b_5^2\ge 1)\,.
$$
Setting
\begin{equation}\label{4.27}
\begin{array}{l}
b_5^2=\{4\sin^2\phi\}/\{\kappa_{22}+\kappa_{33}+(\kappa_{23}+\kappa_{32})\cos\phi\}\\
\mu_2=\kappa_{11}b_5^2
\end{array}
\end{equation}
yields
$$
\P(m(t)b_5^2\ge 1)=\P(J_4(W_1,W_2,W_3;t)\ge 1)\,,
$$
where $J_4$ is defined in (\ref{2.18}), that is, (\ref{2.21}) holds for $\theta\in\Gamma_5(\rho)$.

One can check that for $\theta\in\Gamma_6(\rho)$ the limiting distribution of $\tau(h)$ coincides
 with that for $\theta\in\Gamma_5(\rho)$ with $b_6=b_5$.
 
By a similar argument it can be shown that for $\theta\in\Gamma_7(\rho)$
$$
\lim_{h\to\infty}\P(\tau(h)\le b_7t\sqrt{h})=\P(\tilde{m}(t)b_7^2\ge 1)\,,
$$ 
where
$$
\tilde{m}(t)=\tilde{\kappa}_{11}J_1(W_3;t)+\tilde{\kappa}_{22}J_1(W_4;t)+
$$
$$
+\frac{J_2(W_1,W_2;t)}{4\sin^2\phi}\left(\tilde{\kappa}_{33}
+\tilde{\kappa}_{44}+(\tilde{\kappa}_{34}+\tilde{\kappa}_{43})\cos\phi\right)\,,
$$
$\tilde{\kappa}_{ij}=\sigma^2<(Q_7Q_7^{\prime})^{-1}>_{ij},\ Q_7$ is corresponding transformation matrix which relates the original process with the decomposed ones.

Setting
\begin{equation}\label{4.28}
\begin{array}{l}
b_7^2=(4\sin^2\phi)/\left(\tilde{\kappa}_{33}+\tilde{\kappa}_{44}+(\tilde{\kappa}_{34}+\tilde{\kappa}_{43})\cos\phi\right)\,,\\
\mu_3=\tilde{\kappa}_{11}b_7^2,\ \mu_4=\tilde{\kappa}_{22}b_7^2
\end{array}
\end{equation}
one comes to (\ref{2.21}). This completes the proof of Theorem 2.2. $\endproof$



\begin{flushright}
\begin{tabular}{lcl}
   L.Galtchouk                       &$\quad$& V.Konev              \\
 IRMA, Department of Mathematics     &$\quad$& Department of Applied Mathematics and\\      
 Strasbourg University               &$\quad$& Cybernetics, Tomsk University        \\
 7, str. Rene Descartes              &$\quad$& 36, str. Lenin,             \\
 67084, Strasbourg Cedex             &$\quad$& 634050, Tomsk,\\
 France                              &$\quad$&     Russia\\
 e-mail: galtchou@math.u-strasbg.fr  &$\quad$& vvkonev@vmm.tsu.ru         \\
\end{tabular}
\end{flushright} 

\end{document}